\documentclass[envcountsame,orivec,runningheads,11pt]{llncs}
\usepackage[utf8]{inputenc}
\usepackage{a4}
\usepackage{a4wide,amsmath, amssymb, amsfonts,xspace,stmaryrd}
\usepackage{hyperref}
\usepackage{graphicx}
\usepackage{footmisc}
\makeatletter
\DeclareRobustCommand*{\bfseries}{%
  \not@math@alphabet\bfseries\mathbf
  \fontseries\bfdefault\selectfont
  \boldmath
}
\makeatother

\newtheorem{myexample}[theorem]{Example}
\newtheorem{fact}[theorem]{Fact}
\newtheorem{myremark}[theorem]{Remark}

\newcommand{\COMMENTED}[1]{}

\newcommand{\moc}{\operatorname{moc}}
\newcommand{\dom}{\operatorname{dom}}

\newcommand{\myrho}{\rho}
\newcommand{\myrhob}{\myrho_{\textup{b}}}
\newcommand{\myrhol}{\myrho_<}
\newcommand{\myrhog}{\myrho_>}

\newcommand{\calC}{\mathcal{C}}

\newcommand{\ID}{\mathbb{D}}
\newcommand{\IR}{\mathbb{R}}
\newcommand{\IP}{\mathbb{P}}

\newcommand{\IN}{\mathbb{N}}

\newcommand{\sdzero}{\textup{\texttt{0}}\xspace}
\newcommand{\sdone}{\textup{\texttt{1}}\xspace}

\newcommand{\calA}{\mathcal{A}}
\newcommand{\calB}{\mathcal{B}}
\newcommand{\Cantor}{\mathcal{C}}

\newcommand{\cadlag}{c\`{a}dl\`{a}g\xspace}
\newcommand{\caglad}{c\`{a}gl\`{a}d\xspace}

\usepackage[usenames]{color}  

\makeatletter
\newcommand*\bigdot{\mathpalette\bigdot@{.5}}
\newcommand*\bigdot@[2]{\mathbin{\vcenter{\hbox{\scalebox{#2}{$\m@th#1\bullet$\,\,}}}}}
\makeatother

\title{Randomized Computation of Continuous Data: \\ Is Brownian Motion Computable?\thanks{%
Supported by the National Research Foundation of Korea 
(grant NRF-2017R1E1A1A03071032),
by the International Research \& Development Program of
the Korean Ministry of Science and ICT (grant NRF-2016K1A3A7A03950702),
and by the \protect\raisebox{-1pt}{\protect\includegraphics[height=8pt]{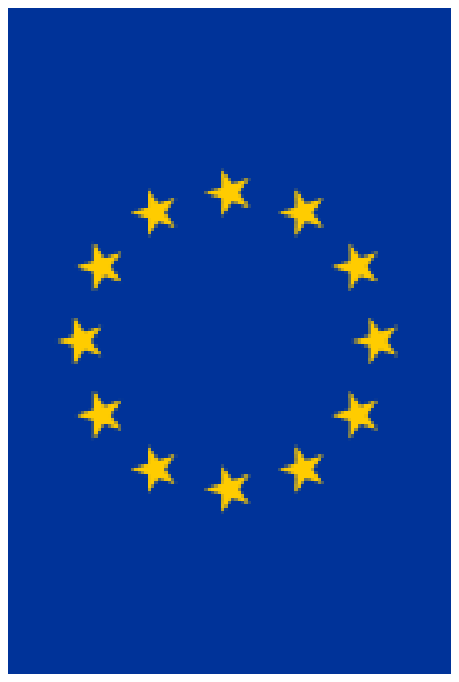}}~European Union's
Horizon 2020 MSCA IRSES project 731143. 
Dedicated to the memory of Klaus Keimel who in 2014 suggested to the last author to study the computability of \cadlag functions.
We thank Volker Betz and Frank Aurzada for advice and assistance.}}
\titlerunning{Randomized Computation of Continuous Data: Is Brownian Motion Computable?}
\authorrunning{Willem Fouch\'{e}, Hyunwoo Lee, Donghyun Lim, Sewon Park, Matthias Schr\"oder, Martin Ziegler}
\author{Willem L. Fouch\'{e}$^1$, Hyunwoo Lee$^2$, Donghyun Lim$^2$, \\
Sewon Park$^2$, Matthias Schr\"oder$^3$, Martin Ziegler$^2$}
\institute{$^1$ University of South Africa \qquad $^2$ KAIST \qquad $^3$ University of Birmingham}
\date{\textbf{keywords:} Computer Science of Continuous Data, Type-2 Theory of Effectivity, Randomization, Brownian Motion}

\begin{document}
\maketitle

\begin{abstract}
We consider randomized computation of continuous data in the sense of Computable Analysis.
Our first contribution formally confirms that it is no loss of generality to 
take as sample space the Cantor space of infinite \emph{fair} coin flips.
This extends [Schr\"{o}der\&Simpson'05] and [Hoyrup\&Rojas'09] considering 
sequences of suitably and adaptively \emph{biased} coins.

Our second contribution is concerned with 1D \emph{Brownian Motion} (aka Wiener Process),
a probability distribution on the space of continuous functions $f:[0,1]\to\IR$ with $f(0)=0$
whose computability has been conjectured [Davie\&Fouch{\'e}'2013; arXiv:1409.4667,\S6].
We establish that this (higher-type) random variable is computable iff
some/every computable family of moduli of continuity (as ordinary random variables)
has a computable probability distribution with respect to the Wiener Measure.
\end{abstract}


\section{Introduction}
Randomization is a powerful technique in classical (i.e. discrete) Computer Science:
Many difficult problems have turned out to admit simple solutions by algorithms that `roll dice'
and are efficient/correct/optimal with high probability
\cite{DBLP:journals/siamcomp/DietzfelbingerKMHRT94,Berenbrink1999,DBLP:journals/rsa/CzumajS00,DBLP:journals/jcss/BeierV04}.
Indeed, fair coin flips have been shown computationally universal \cite{Walker:1977:EMG:355744.355749}.
Over continuous data, well-known closely connected to topology \cite{Grz57} \cite[\S2.2+\S3]{Wei00}, 
notions of probabilistic computation are more subtle \cite{BGH15a,CollinsArxiv2015}.

\subsection{Overview}
Section~\ref{s:Represent} resumes from \cite{SS06} the question of how to represent Borel probability measures.
\cite[Proposition~13]{SS06} had established that, on `reasonable' spaces, 
every such distribution can be represented by the distribution of an infinite sequence
of coin flips (i.e. over Cantor space) with a \emph{suitably} and adaptively biased coin.
Theorem~\ref{t:Universal} shows that such can in turn be represented by `fair' coins.
Lemma~\ref{l:Computable} characterizes computability of a Borel probability measure on the reals:
Necessary and sufficient is that both the lower and upper semi-inverse of its cumulative
probability distribution are, respectively, lower and upper semi-computable real functions.

Section~\ref{s:Brownian} approaches the question of whether Brownian Motion
(aka Wiener Process), a popular probability distribution on the space of
continuous real functions, is computable:
in the strong sense of Subsection~\ref{ss:Computability} underlying \cite{DF13,MTY13}.
Subsection~\ref{ss:Hyunwoo} recalls several known mathematical characterizations
of this distribution, and relates their types of convergence to 
weaker notions of probabilistic computation \cite{Bos08} while
pointing out their differences to the strong sense.
It turns out that quantitative continuity of Brownian Motion, 
captured in terms of some modulus considered as a derived random variable, 
constitutes the major obstacle:
Theorem~\ref{t:Brownian} establishes that computability of the 
probability distribution of any computable such a modulus is both sufficient
and necessary for the computability of Brownian Motion.
This reduces the conjecture from the probability distribution on a function
space to that of an ordinary real random variable.

\section{Representing Borel Probability Measures}
\label{s:Represent}
Recall that a \emph{measure space} is a triple $(X,\calA,\mu)$, where $X$ is a non-empty set,
$\calA$ is a $\sigma$-algebra over $X$, and $\mu$ is a measure on $(X,\calA)$.
For measure spaces $(X,\calA,\mu)$ and $(Y,\calB,\nu)$ and a measurable partial mapping $F:\subseteq X\to Y$,
$\nu$ is the \emph{pushforward} measure of $\mu$ w.r.t. $F$ 
if $\mu\big(F^{-1}[V]\big)$ is defined and equal to $\nu(V)$ for every $V\in\calB$.
In this case we say $F$ \emph{realizes} $\nu$ on $\mu$ 
and write $\nu\preccurlyeq\mu$. 
This notion is similar to, but not in danger of confusion
with, \cite[Definition~2.3.2]{Wei00}; we will generalize
it in Definition~\ref{d:Computable}.
Note that realizability is transitive; 
and a realizer $F$ must have $\dom(F)\in\calA$ of measure $\nu(Y)$.

The Type-2 Theory of Effectivity 
employs Cantor space to encode, and define computation over, 
any topological T$_0$ space, such as real numbers 
and continuous real functions \cite[\S3.2+\S4]{Wei00}.

\begin{figure}[htb]
\begin{center}
\includegraphics[width=0.9\textwidth]{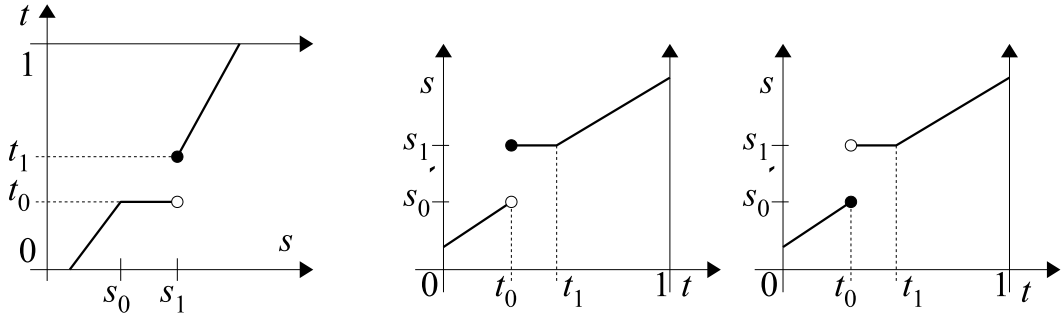}
\caption{\label{f:cadlag}Example cumulative distribution function with upper/lower semi-inverse}
\end{center}
\end{figure}

\begin{myexample}
\label{x:Probability}
\begin{enumerate}
\item[a)]
Consider the real unit interval $X=[0,1]$ equipped with the $\sigma$-algebra $\calA$ of Borel subsets
and the Lebesgues probability measure $\lambda$.
\item[b)]
Consider Cantor space $\Cantor=\{\sdzero,\sdone\}^\omega$ equipped with the $\sigma$-algebra $\calB$ of Borel subsets
and the canonical (=fair coin flip) probability measure $\gamma$:
$\gamma(\vec w\circ\Cantor)=2^{-|\vec w|}$, where $|\vec w|=n$ denotes
the length of $\vec w=(w_0,\ldots,w_{n-1})\in\{\sdzero,\sdone\}^n$.
\item[c)]
The continuous total mapping $\myrhob:\Cantor\ni\bar b\mapsto \sum_{j\geq0} b_j2^{-j-1}\in[0,1]$ 
realizes $([0,1],\calA,\lambda)$ on $(\Cantor,\calB,\gamma)$:
$([0,1],\calA,\lambda)\preccurlyeq(\Cantor,\calB,\gamma)$.
\item[d)]
Consider the real line $\IR$ equipped with (the Borel $\sigma$-algebra and)
the standard Gaussian/normal probability distribution,
realized on $\lambda$ via the partial mapping
$G:(0,1)\ni t\mapsto \Phi^{-1}(t)\in\IR$ for the cumulative distribution 
\[ \Phi\;:\;\IR\;\ni\; s \;\mapsto \int\nolimits_{-\infty}^s \exp(-t^2/2) /\sqrt{2\pi}\, dt \;\in\; [0,1] \]
\item[e)]
Consider $[0,1]$ equipped with the Dirac point measue $\delta_r$ for some $r\in(0,1)$.
It is realized on $\big([0,1],\calA,\lambda\big)$ via 
the constant function $H:[0,1]\mapsto \{r\}$.
\item[f)]
Consider $[0,1]$ equipped with the Cantor measure.
It is realized on $\big([0,1],\calA,\lambda\big)$ via 
the inverse of \emph{Devil's Staircase} (aka Cantor–Vitali function).
\end{enumerate}
\end{myexample}
Let us combine and generalize the real
Items~(d), (e), and (f) of Example~\ref{x:Probability}:

\begin{lemma}[Real Case]
\label{l:Cadlag}
Fix $X=\IR$ equipped with some Borel probability measure $\mu$.
Recall that its cumulative distribution function
$\IR\ni s\mapsto \mu\big((-\infty,s]\big)\in[0,1]$
is \emph{\cadlag}
(continuous from right with left limits)
and non-decreasing, hence upper semi-continuous.
On the other hand
$s\mapsto \mu\big((-\infty,s)\big)$ is \emph{\caglad}
and lower semi-continuous.
Now consider the cumulative distribution function's
\emph{upper} and \emph{lower semi-inverse}:
\begin{align*}
F^\mu_>: (0,1)\;\ni\; t\;\mapsto\; &
\inf\big\{ s \in\IR \:\big|\: \mu\big((-\infty,s)\big)>t\big\}
\;=\;
\min\big\{ s \:\big|\: \mu\big((-\infty,s]\big)>t\big\}  \\
=\; &
\sup\big\{ s \in\IR \:\big|\: \mu\big((-\infty,s]\big)\leq t\big\} 
\\[0.5ex]
F^\mu_<: (0,1)\;\ni\; t\;\mapsto\; &
\max\big\{ s \in\IR \:\big|\: \mu\big((-\infty,s)\big)<t\big\}
\;=\;
\sup\big\{ s \:\big|\: \mu\big((-\infty,s]\big)<t\big\} \\
=\; &
\sup\big\{ s \:\big|\: \mu\big((s,\infty)\big)> 1-t\big\} 
\end{align*}
$F^\mu_>$ is \cadlag and upper-semicontinuous;
$F^\mu_<$ is \caglad and lower-semicontinuous;
and both realize $\mu$ on $\big([0,1],\calA,\lambda\big)$:
see Figure~\ref{f:cadlag}.
\end{lemma}
\noindent
In 
the sequel we consider topological spaces,
implicitly equipped with the Borel $\sigma$-algebra,
and a Borel probability measure.
\cite[Proposition~13]{SS06} establishes the following:

\begin{fact}
\label{f:SchroederSimpson06}
To every 2nd countable T$_0$ space $X$ with Borel probability measure $\mu$
there exists a Borel probability measure $\tilde\gamma$ on $\Cantor$
such that $(X,\mu)$ has a continuous partial realizer over $(\Cantor,\tilde\gamma)$.
\end{fact}
The metric case is treated in \cite[Theorem~5.1.1]{HR09}.
We show that the probability measure $\tilde\gamma$ on $\Cantor$ can 
in fact be chosen as the canonical `fair' one:

\begin{theorem}
\label{t:Universal}
Every Borel probability measure $\tilde\gamma$ on Cantor space $\Cantor$
admits a continuous partial realizer over the `fair' measure $(\Cantor,\gamma)$.
The realizer is defined on $\Cantor$ with the exception of at most countably many points.
\end{theorem}
Indeed, Fact~\ref{f:SchroederSimpson06} and transitivity together imply 
that every 2nd countable T$_0$ space with a Borel probability measure
to admit a continuous partial realizer over $(\Cantor,\gamma)$.
One cannot hope for a total realizer in general, though:

\begin{proposition}
\label{p:Schroeder}
Let $\mu$ be a Borel probability measure on $\Cantor$ such that there
is some $\bar v \in \Cantor$ such that the measure of the basic open set
$\bar v\,\Cantor$ is \emph{non-}dyadic.
Then there is no total continuous function $F: \Cantor \to\Cantor$ with $\gamma\circ F^{-1}=\mu$.
\end{proposition}
%

\subsection{Proofs of Theorem~\ref{t:Universal} and Proposition~\ref{p:Schroeder}}

\begin{proof}[Theorem~\ref{t:Universal}]
For each open interval $I=(a,b)\subseteq[0,1]$
consider the set $\Cantor_I=\myrhob^{-1}[I]\subseteq\Cantor$ 
of measure $\gamma(\Cantor_I)=\lambda(I)$.
Note that $\Cantor_{I\cup J}=\Cantor_I\cup\Cantor_J$ 
and $\Cantor_{I\cap J}=\Cantor_I\cap\Cantor_J$.
Fix $n\in\IN$ and equip $\{\sdzero,\sdone\}^n$ with the total lexicographical order; 
and consider the disjoint open intervals 
\[ I_{\vec 0}\:=\:\big(0,\tilde\gamma(\vec 0\circ\Cantor)\big)
\quad\text{as well as}\quad
I_{\vec w}\:=\:\Big(\sum\nolimits_{\vec v\pmb{<}\vec w} \tilde\gamma(\vec v\circ\Cantor),
\sum\nolimits_{\vec v\pmb{\leq}\vec w} \tilde\gamma(\vec v\circ\Cantor)\Big) \]
of lengths $\lambda(I_{\vec w})=\tilde\gamma(\vec w\circ\Cantor)$
for each $\vec w\in\{\sdzero,\sdone\}^n\setminus\vec 0$.
Since $\tilde\gamma$ is a Borel probability measure on $\Cantor$, 
these lengths add up to $\sum_{\vec w} \tilde\gamma(\vec w\circ\Cantor)=1$.
Also note that $I_{\vec w\,\sdzero},I_{\vec w\,\sdone}\subseteq I_{\vec w}$
are disjoint with lengths $\lambda(I_{\vec w\,\sdzero})+\lambda(I_{\vec w\,\sdone})=\lambda(I_{\vec w})$;
and that $I_{\vec w}$ may be empty in case $\tilde\gamma(\vec w\circ\Cantor)=0$.
Finally abbreviate 
\[ \Cantor_{\vec w}:=\Cantor_{I_{\vec w}} \quad\text{ and }\quad
F_n:\subseteq\Cantor\to\{\sdzero,\sdone\}^n, \quad F_n\big|_{\Cantor_{\vec w}}:\equiv\vec w \]
so that $F_n$ is defined except for at finitely many arguments
(namely the binary encodings of the real interval endpoints)
with $F_n^{-1}(\vec w)=\Cantor_{\vec w}$
of measure $\gamma\big(F_n^{-1}(\vec w)\big)=\tilde\gamma\big(\Cantor_{\vec w}\big)$.
Since $F_{n+1}(\bar u)\in F_n(\bar u)\circ\{\sdzero,\sdone\}$,
$F(\bar u):=\lim_n F_n(\bar u)\in\Cantor$ is well-defined 
(except for at countably many arguments) and continuous
with $F^{-1}[\vec w\,\Cantor]=\Cantor_{\vec w}$ for every $\vec w\in\{\sdzero,\sdone\}^*$.
Hence $\gamma\circ F^{-1}$ coincides with $\tilde\gamma$ on the basic clopen subsets of $\Cantor$
and, being Borel measures, also on all Borel subsets.
\qed\end{proof}

\begin{proof}[Proposition~\ref{p:Schroeder}]
Suppose that $F$ is such a continuous function. For every $n$ there is a word
$\vec w_n$ of length $n$ such that $F[\vec w_n\,\Cantor]$ contains an element of
$\vec v\,\Cantor$ and an element of its complement, because otherwise the
preimage $\bar v\,\Cantor$ would be the finite union of all open balls
$\vec w\,\Cantor$ with all $\vec w$ of length $n$ satisfying $F[\vec w\,\Cantor]
\subseteq \vec v\,\Cantor$; but the $\gamma$-measure of this union is dyadic.
By the fan theorem (or by the fact the $\Cantor$ is a (sequentially)
compact space), there is some $\bar p \in\Cantor$ and some infinite subset $I$
of $\IN$ such that $\vec w_i$ is a prefix of $\bar p$ for all $i \in I$. But $F$ cannot
be continuous in the point $\bar p$, because no prefix of $\bar p$ can tell whether
$F(\bar p)$ is inside or outside the clopen set $\vec v\,\Cantor$, a contradiction!
\qed\end{proof}

\subsection{Computability of Borel Probability Distributions}
\label{ss:Computability}

Of course a realizer in the sense of Theorem~\ref{t:Universal} is usually far from unique.
We are interested in those computable with respect to a representation of the space under consideration:
\begin{definition}
\label{d:Computable}
Fix a Borel probability measure $\mu$ on $X$
and a representation $\xi:\subseteq\Cantor\twoheadrightarrow X$
in the sense of TTE {\rm\cite[\S3]{Wei00}}.
A \emph{$\xi$-realizer} of $\mu$ is a mapping $G:\subseteq\Cantor\to\dom(\xi)$
such that $\xi\circ G:\subseteq\Cantor\to X$ is a realizer of $\mu$
(over the `fair' measure) in the above sense.
Call $\mu$ \emph{$\xi$-computable} if it has a computable $\xi$-realizer.
\end{definition}
Note that $\dom(\xi)\subseteq\Cantor$ is a 2nd countable T$_0$ space,
equipped with the pushforward measure of $\xi$:
hence Fact~\ref{f:SchroederSimpson06} and Theorem~\ref{t:Universal}
together assert that a (possibly uncomputable) $\xi$-realizer exists!

\begin{myexample}
\label{x:Computable}
Recall {\rm\cite[Definition~2.3.2]{Wei00}} computable reduction between representations
as well as the importance of \emph{admissible} ones {\rm\cite[\S3.2]{Wei00}},
which $\myrhob$ does not belong to {\rm\cite[Theorem~4.1.13.6]{Wei00}}.
\begin{enumerate}
\item[a)]
The identity on $\dom(\myrhob)\subseteq\Cantor$ constitutes
a computable $\myrhob$-realizer of the Lebesgues measure on $[0,1]$
according to Example~\ref{x:Probability}c).
\item[b)]
If $\mu$ is $\xi$-computable and if $\xi\preccurlyeq\xi'$ holds,
then $\mu$ is also $\xi'$-computable.
\item[c)]
In particular the Lebesgues measure on $[0,1]$
is $\myrho$-computable for the 
admissible representation 
$\myrho$ {\rm\cite[Theorem~4.1.13.7]{Wei00}}.
\item[d)]
And so is the Gaussian distribution from
Example~\ref{x:Probability}d)
as well as the Cantor distribution from
Example~\ref{x:Probability}f).
The Dirac distribution $\delta_r$ is $\myrho$-computable
~iff~ $r$ is $\myrho$-computable.
\item[e)]
For admissible $\xi:\subseteq\Cantor\twoheadrightarrow X$,
integration $\calC(X)\ni f\mapsto\int_X f(x)\,d\mu(x)\in\IR$
is $\big([\xi\!\to\myrho],\myrhol\big)$-computable 
iff $\mu$ has a computable $\xi$-realizer {\rm\cite{Sch07b}}.
\end{enumerate}
\end{myexample}
Regarding (eq), recall {\rm\cite[\S3.3]{Wei00}} that every admissible representation
$\xi$ of $X$ induces a canonical admissible representation $[\xi\!\to\myrho]$
of the space $\calC(X)$ of continuous functions $f:X\to\IR$;
and recall {\rm\cite[Lemma~4.1.8]{Wei00}} the representation $\myrhol$ of $\IR$ 
encoding approximations from below only. 
It is well-known that every $(\myrho,\myrhol)$-computable function
must be lower semi-continuous \cite{WZ00,Zie07a}.
Resuming Lemma~\ref{l:Cadlag},
we can now characterize the real case \cite{Wei99a}:

\begin{lemma}
\label{l:Computable}
Fix a Borel probability measure $\mu$ on $\IR$
with cumulative distribution function and 
lower and upper semi-inverse $F^\mu_<$ and $F^\mu_>$.
Then $\mu$ is computable in the sense of Definition~\ref{d:Computable}
~iff~ both $F^\mu_<$ is $(\myrho,\myrhol)$-computable
and $F^\mu_>$ is $(\myrho,\myrhog)$-computable.
\end{lemma}

\section{Characterizing Computability of Brownian Motion}
\label{s:Brownian}

1D \emph{Brownian Motion} aka \emph{Wiener Process}
is a probability measure on the space $X:=\calC_0[0,1]$ 
of (i) continuous functions $W:[0,1]\to\IR$ satisfying (ii) $W(0)=0$
and characterized by the following properties:
\begin{enumerate}
\item[iii)]
For every $0\leq r<s<t\leq 1$,
$W(t)-W(s)$ is independent of $W(r)$.
\item[iv)]
$W(t)$ is Gaussian normally distributed with mean $W(s)$
and variance $|t-s|$.
\end{enumerate}
Compare \cite{LeGall,KarlinTaylor,Maler2015} for details.
Here we approach the question of whether this measure is computable \cite{DBLP:journals/tcs/Fouche08,DF13}
in the sense of admitting a \emph{computable} $[\myrho\!\to\!\myrho]$-realizer;
recall Definition~\ref{d:Computable}.

\begin{figure}[htb]
\begin{center}
\includegraphics[width=0.9\textwidth,height=0.2\textheight]{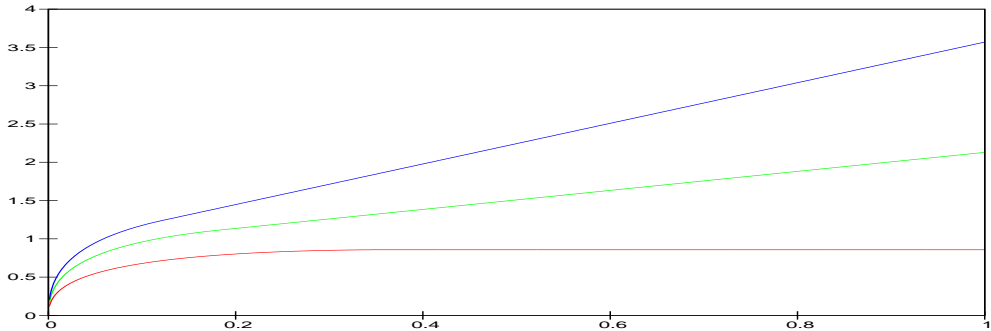}
\caption{\label{f:modulus}Parameterized moduli of continuity of an increasing sequence of subsets of Brownian Motion}
\end{center}
\end{figure}

\begin{myremark}
\label{r:Continuity}
The representation $[\myrho\!\to\!\myrho]$ encodes any $f\in\calC[0,1]$ 
via both (I) its values $f(a/2^n)$ on the countable dense subset
of dyadic rationals $\ID:=\bigcup_n\ID_n$, $\ID_n:=\big\{0/2^n,1/2^n,\ldots,2^n/2^n\big\}$,
and (II) a binary modulus of continuity of $f$:
a sequence $\moc:\IN\to\IN$ such that $|s-t|\leq2^{-\moc(n)}$ implies $|W(s)-W(t)|\leq2^{-n}$;
see {\rm\cite[\S6.1]{Wei00}}.
\begin{enumerate}
\item[a)]
Based on Example~\ref{x:Probability}e),
Conditions~(iii) and (iv) immediately yield an algorithm
for computably `guessing' the values $W\big|_{\ID}$ 
according to (I) iteratively on-the-fly 
with respect to the appropriate Gaussian normal distribution
in relation to the previous values.
However this approach does not allow to then 
(II) determine $\moc(n)$ within finite time:
with small but positive probability, $W\big|_{\ID_{n+m}}$
may exceed any purported upper bound $\moc(n)$.
\item[b)]
Conversely first (II) `guessing' $\moc(n)$ 
requires to know the probability distribution of the random variable $\moc$
\emph{exactly}: otherwise the resulting Wiener Process will have skewed quantitative continuity.
This in turn affects (I) the distribution of $W\big|_{\ID}$,
with properties (iii) and (iv) now having probabilities \emph{conditional} to said $\moc$. 
Recall the following generic (though not necessarily efficient) way of modifying any 
randomized algorithm to adjust its internal guesses to become conditional to some event $E$:
For every independent sample $s$, test shether $g\in E$; if not, discard $s$ and sample again
--- until obtaining one that complies with $E$.
\item[c)]
By \emph{L\'{e}vy's modulus of continuity theorem}, 
with probability 1 it holds
\begin{equation} 
\label{e:Hoelder}
\lim_{h\to0}  \: \sup_{|s-t|\leq h} \: \frac{|W(s)-W(t)|}{\sqrt{2h\ln(1/h)}} \;=\; 1 
\end{equation}
The Wiener Process is thus $\alpha$-H\"{o}lder continuous for every exponent $\alpha>1/2$, but not for $\alpha=1/2$.
\item[d)]
More explicitly,
abbreviating $E:=\exp(1)$ and $y_c:=\sqrt{2\ln(Ec)/c}$, 
Equation~\eqref{e:Hoelder} says that, to every $W$ (except for a subset of measure zero)
there exists some least $c=c(W)\geq1$ such that
\begin{equation}
\label{e:Levy}
\omega(h,c)\;:=\;
\left\{ \begin{array}{cl} 
\sqrt{2ch\ln(1/h)} &:\; h\leq 1/Ec \\
y_c+(h-1/Ec)\cdot c\cdot\ln(c)/y_c  &:\; h\geq 1/Ec  
\end{array}\right.
\end{equation}
depicted in Figure~\ref{f:modulus}
constitutes a parameterized modulus of continuity of $W$ in the following sense:
\item[e)]
For a function $f:X\to Y$ between metric spaces $(X,d)$ of diameter 1 and $(Y,e)$,
a (classical, as opposed to binary) modulus of continuity is a mapping $\omega:[0,1]\to[0,\infty)$
such that it holds
\begin{equation}
\label{e:Modulus}
\omega(0)=0 \quad\wedge\quad \forall x,x'\in X: \quad e\big(f(x),f(x')\big) \;\leq\; \omega\big(d(x,x')\big)
\enspace .
\end{equation}
If $f$ is continuous with compact domain, then it has a modulus of continuity $\omega$.
It $X$ is additionally convex, $\omega$ can be chosen subadditive.
\end{enumerate}
\noindent The $c\geq1$ from Item~(c) is thus an unbounded real random variable,
parameterizing the family of subadditive moduli of continuity from Equation~\eqref{e:Levy}
strictly increasing in both arguments.
\end{myremark}
We can now state our main result characterizing computability of the Wiener Process
in terms of computability of the probability distribution of any/all 
parameterized moduli of continuity:

\begin{theorem}
\label{t:Brownian}
Let $\omega:[0,1]\times[1,\infty)\to[0,\infty)$ denote any computable (and thus continuous) 
one-parameter family of subadditive functions strictly increasing in both arguments with $\omega(0,c)\equiv0$.
Suppose that to every Wiener Process $W$ (except for a subset of measure zero)
there exists a (necessarily unique) least $c=c(W)\geq1$ such that $\omega(\:\cdot\:,c)$ 
constitutes a modulus of continuity of $W$ in the sense of Remark~\ref{r:Continuity}e).
Then the following are equivalent:
\begin{itemize}
\item[\textbullet] The Wiener Process $W$ is computable (formally: has a computable $[\myrho\!\to\!\myrho]$-realizer).
\item[\textbullet] The random variable $c$ has a computable probability distribution.
\item[\textbullet] There exists a random variable $\tilde c$ with computable probability distribution
such that $\omega(\:\cdot\:,\tilde c)$ is a modulus of continuity of $W$ with probability 1.
\end{itemize}
\end{theorem}
%

\subsection{Na\"\i{}ve Approaches and their Deficiencies}
\label{ss:Hyunwoo}
Before proceeding to the proof of Theorem~\ref{t:Brownian},
let us report some well-known alternative characterizations
of mathematical Brownian Motion
and why they do not imply computability.

\begin{myexample}
\label{x:Wiener}
For probability spaces $(X,\mu)$ and $(Y,\nu)$,
recall that a sequence $R_n:\subseteq Y\to X$ of random variables 
converges \emph{almost surely} to $R:\subseteq Y\to X$ if the set
$\big\{y : R_n(y)\to R(y)\big\}\subseteq Y$ has $\nu$-measure 1.

On the other hand for $(X,d)$ a metric space,
\emph{uniform} almost sure convergence of $R_n$ to $R$
means that there exists $U\subseteq\dom(R)\cap\bigcap_n\dom(R_n)$ of $\nu$-measure 1
such that $\sup_{y\in U} d\big(R_n(y),R(y)\big)\to0$.

\begin{enumerate}
    \item [a)] Let $\varphi_0(t)=t$ and
\[
\varphi_{n,j}(t)=
\begin{cases}
2^{(n-1)/2}\cdot(t-\frac{k-1}{2^n})&   \frac{k-1}{2^n}\leq t \leq \frac{k}{2^n} \\
2^{(n-1)/2}\cdot(\frac{k+1}{2^n}-t)&   \frac{k}{2^n}\leq t \leq \frac{k+1}{2^n} \\
0& \text{otherwise}
\end{cases}, \qquad 0\leq k<2j, \quad 1\leq j\leq 2^{n-1}
\]
denote the \emph{Schauder} `hat' functions and $R_{n,j}$ independent standard normally distributed random variables. Then following sequence converges to the Wiener Process almost surely:
\begin{gather}
W(t)^N(\omega)=R_{0}(\omega)t+\sum_{n=1}^{N}\sum_{j=1}^{2^{n-1}}R_{n,j}(\omega)\varphi_{n,j}(t)
\end{gather}

\item[b)]  Let $R_i$ be independent standard normally distributed random variables
(Example~\ref{x:Probability}d). Then following sequence converges to the Wiener Process in mean.

\begin{gather}
W(t)^N(\omega)=\sqrt{2} \sum_{i=1}^{N} R_i \frac{\sin{(k-\frac{1}{2})\pi t}}{(k-\frac{1}{2}) \pi}
\end{gather}

\item[c)] Let $(X_i)_{i\in \IN}$ be independent random variable with mean 0 and variance 1 and $S_n=\sum_{i=1}^{n} X_i$. 
Then following sequence converges to the Wiener Process in distribution:
\begin{gather}
    W(t)^N(\omega)=\frac{S_{\lfloor nt \rfloor }}{\sqrt{N}}
\end{gather}
\end{enumerate}
\end{myexample}

\subsection{Proof of Theorem~\ref{t:Brownian}}
\label{ss:Proof}
We first record that the hypotheses 
ensure that $\omega(\:\cdot\:,C)$ does constitute a modulus of continuity
for every $C\geq1$, namely one of $\omega(\:\cdot\:,C)$ itself.
Moreover strict monotonicity in $C$ asserts that the
measure of all \emph{those} Wiener Processes $W$ which have $\omega(\:\cdot\:,C)$
as modulus of continuity is strictly increasing and continuous.
Hence we can apply Lemma~\ref{l:Computable}
with continuous (as opposed to just \cadlag) cumulative probability distribution
and with lower and upper semi-inverse coinciding and continuous.

\bigskip\noindent
First suppose the random variable $\tilde c$ parameterizing $\omega$
has a computable inverse cumulative probability distribution.
Similarly to Example~\ref{x:Probability}d),
this allows to algorithmically `guess' the value $\tilde C$ of $\tilde c$ according to said distribution;
and computability of $\omega$ can be turned into an (upper bound on the) 
binary modulus of continuity $\moc$ (II). 
Regarding (I),
having guessed and fixed a modulus of continuity $\omega(\:\cdot\:,\tilde C)$
affects properties (iii) and (iv) of the Wiener Process.
As mentioned in Remark~\ref{r:Continuity}b),
this can be atoned for by discarding guesses for values
$W(t)$ that violate $\omega(\:\cdot\:,\tilde C)$ --- but is complicated in our case 
with undecidable real equality \cite[Exercise~4.2.9]{Wei00}.
So we make a point of carefully using only \emph{strict} inequalities,
which are at least semi-decidable:

\begin{quote}
Beginning with $V:=\emptyset$ iteratively/on-demand 
guess a new value $W(s)$, $s\in\ID$, subject to (iii) and (iv).
Then check whether it complies with all previously guessed values $W(t)$, $t\in V\subseteq\ID$,
in satisfying $|W(s)-W(t)|\pmb{<}\omega(|s-t|,\tilde C)$; 
and if so, add $s$ to $V$.
On the other hand if there is some $t\in V$ with 
$|W(s)-W(t)|\pmb{>}\omega(|s-t|,\tilde C)$,
then discard and guess again the value of $W(s)$.
\end{quote}

\noindent
Note that the above comparisons exclude and fail in the case
$|W(s)-W(t)|\pmb{=}\omega(|s-t|,\tilde C)$:
which occurs only with probability 0, though:
The above algorithm thus computes $W$ with probability 1.

\bigskip\noindent
Conversely suppose that Brownian Motion has a computable
$[\myrho\!\to\!\myrho]$-realizer $F:\subseteq\Cantor\to\dom([\myrho\!\to\!\myrho])\subseteq\Cantor$.
For each given $C\in\IR$ and $\bar u\in\dom(F)$ and $W:=\big[\myrho\!\to\!\myrho\big](\bar u)$,
computability of $\omega$ then implies \cite[Corollary~6.2.5]{Wei00} computability of 
\[ 
(C,W) \;\mapsto\; \Psi(C,W) \;:= \max_{0\leq s,t\leq1} \omega(|s-t|,C) - |W(s)-W(t)|
\] 
and of $W\mapsto \min\{C: \Psi(C,W) \geq0\}=\max\{C:\Psi(C,W)\leq0\}=c(W)$
defined for almost all $W$:
since $C\mapsto\Psi(C,W)$ is continuous and strictly increasing 
and unbounded by hypothesis.
Again we consider $c(W)=\inf\{C: \Psi(C,W) >0\}=\sup\{C:\Psi(C,W)<0\}$
instead, in order to avoid undecidable real equality.
Preimages of open sets can be exhausted (only) from inside \cite[Theorem~2.4.5.3]{Wei00},
and their measure (only) from below; cmp. Example~\ref{x:Computable}e).
With $F\circ c$ computable and $F$ a realizer of Brownian Motion,
the sought cumulative probability $\IP \big[ W : c(W) \leq C \big]$
thus coincides with the fair/canonical measure $\gamma$ of
$\big(F\circ c\big)^{-1}\big[[0,C)\big]\subseteq\Cantor$,
and with $1-\gamma\Big(\big(F\circ c\big)^{-1}\big[(C,\infty)\big]\Big)$:
the former yields approximations from below, the latter yield approximations from above,
and together they yield approximations up to any given error \cite[Lemma~4.1.9]{Wei00}.
Note that we do not need $\dom(F)$ to be semi-decidable as it has measure zero anyway.
\qed


\bibliographystyle{alpha}
\bibliography{cca,wiener}
\end{document}